\documentclass[final,twoside,11pt,twocolumn,letterpaper]{article}
\usepackage{macro}
\usepackage{graphicx}
\usepackage{amsmath}
\usepackage{times}
\usepackage{amsfonts,amssymb}
\usepackage{amsbsy,amscd}
\usepackage{rotfloat}

\newtheorem{theorem}{Theorem}[section]

\newtheorem{remark}[theorem]{Remark}

\def\R{\mathbb{R}}

\def\ah{H^1(\Omega)}

\def\ahr{H(\rotscal,\,\Omega)}
\def\divrot0{H_0(\dive,\,\rotscal,\,\Omega)}
\def\divrotb0{H_0(\dive,\,\rotscal,\,\Omega_b)}
\def\el{L^2(\Omega)}

\def \dém{{\sl Preuve. }}

\def\H{{\cal H}}

\def\u{\vec{U}}

\newcommand{\ba}{\begin{eqnarray}}
\newcommand{\ea}{\end{eqnarray}}
\newcommand{\basn}{\begin{eqnarray*}}
\newcommand{\easn}{\end{eqnarray*}}
\newcommand{\dive}{\operatorname {div}}

\newcommand{\rots}{\operatorname{curl}}
\newcommand{\rot}{\vec{\,\operatorname{curl\,}}}
\newcommand{\rotscal}{\operatorname{curl}}
\newcommand{\grad}{\vec{\nabla}}
\newcommand{\Dif}[2]{\frac{\mathrm{D} #1}{\mathrm{D} #2}}
\newcommand{\Diff}[2]{\frac{\mathrm{D}^2 #1}{\mathrm{D} {#2}^2}}
\newcommand{\dif}[2]{\frac{\partial #1}{\partial #2}}

\newcommand{\ddt}[1]{\ensuremath{\frac{{\rm d} #1}
{ {\rm d} t}}}
\newcommand{\ddtt}[1][]{\ensuremath{\frac{{\rm d}^2 #1 }{ {\rm d}t^2}}}
\renewcommand{\vec}[1]{\boldsymbol{#1}}

\newcommand{\demi}{\frac{1}{2}}
\newcommand{\dis}{\displaystyle}

\begin{document}
 
\title{NUMERICAL ANALYSIS OF TIME-DEPENDENT GALBRUN EQUATION IN AN INFINITE DUCT}
\author{\underline{K.\ Berriri}$^{\dag,*}$, A.-S.\ Bonnet-Ben Dhia$^{\dag}$, P.\ Joly$^{\dag}$\\
  \textnormal{$^{\dag}$Laboratory POems, UMR 2706 CNRS-INRIA-ENSTA, Paris, FRANCE\\
$^{*}$Email: Kamel.Berriri@inria.fr}}
\date{}
\maketitle
\pagestyle{empty}
\thispagestyle{empty}
\section*{Abstract}
In this paper we are interested in the mathematical and numerical
analysis of the time-dependent Galbrun equation in a rigid duct. This equation models the acoustic propagation in the presence of a flow \cite{Galbrun:1931}. We propose a regularized variational formulation of the problem, in the subsonic case, suitable for an approximation by Lagrange finite elements, and corresponding absorbing boundary conditions.
\section{Introduction}
\noindent Sound propagation in a flow is a subject of great interest
for numerical analysis. The main application concerns noise
reduction in aeronautics. The understanding of phenomena of interaction between acoustic waves and flows is a crucial feature to find the components that efficiently reduce the sound. 

Most of the mathematical and numerical studies for transient
linearized aeroacoustics are based on linearized Euler equations,
whose unknowns are the perturbations of velocity and
pressure. Nevertheless, Galbrun equations, whose unknown is the
lagrangian perturbation of displacement, are an attractive
alternative to model the phenomenon of acoustic propagation in a
flow. Indeed, these equations seem to have a structure similar  to
equations in electromagnetism and elastodynamics and they allow a simple treatment of boundary conditions.

Recent studies on Galbrun equations have focused on time harmonic dependence. We are interested here in studying the transient case.
 
\begin{figure}[H]
\label{image1}
  \centering
  \includegraphics[scale=0.3]{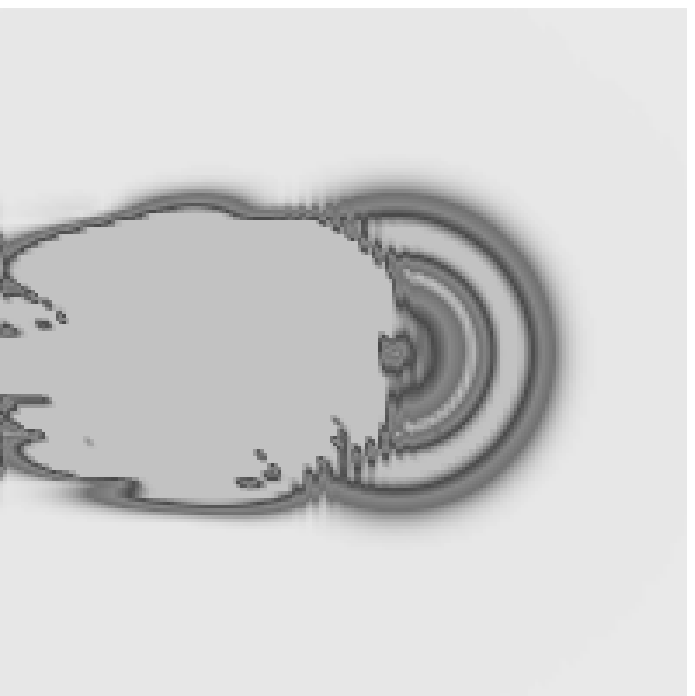} \hspace{0.5cm}
  \includegraphics[scale=0.3]{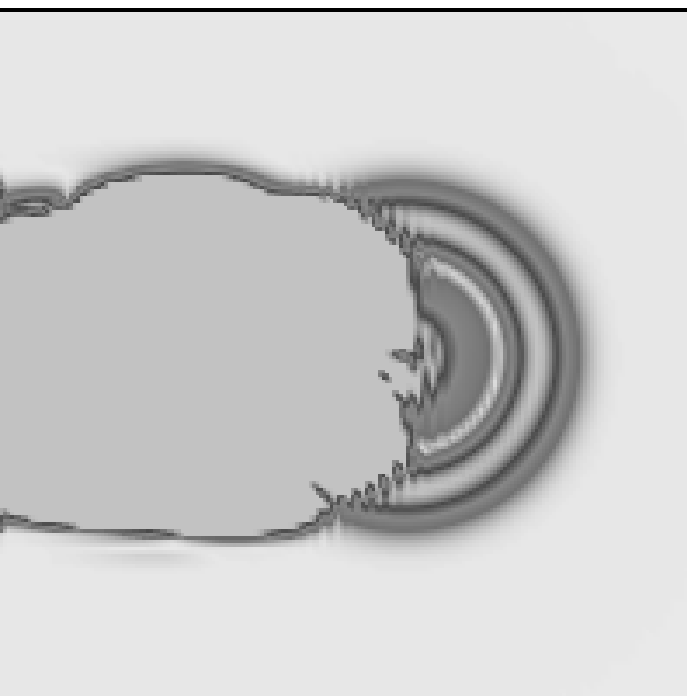} \hspace{0.5cm}
  \includegraphics[scale=0.3]{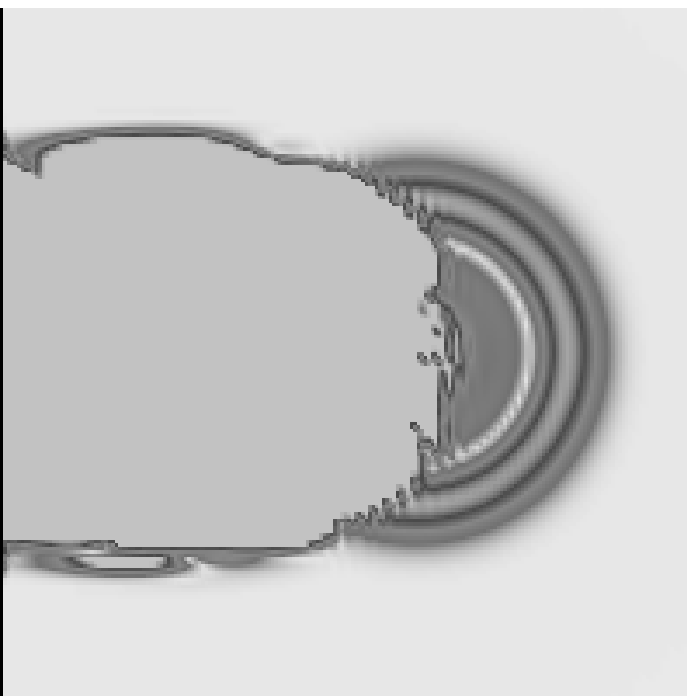}
 \caption{Euclidean norm of $\vec{\xi}$ at  at $t_1=1.5$s and $t_2=1.75$s and $t_3=2$s. Case $s=0$}
  \label{fig2}
\end{figure}

This problem raises many theoretical and numerical difficulties. The
major one is the lack of a natural functional frame for a
variational setting of the problem. Moreover, a naive numerical
resolution of the problem using standard nodal finite elements for the
space discretization is unstable (see Figure \ref{fig2}). In order to overcome
these difficulties , we propose a regularization method analogous to the one
developed for Maxwell's equations \cite{Costabel:1991}. This allows to derive a variational formulation suitable for an approximation by Lagrange finite elements.
\section{Statement of the problem}
\noindent Let $\Omega=\R\times]-h,h[$ be a two-dimensional rigid duct containing a compressible fluid. We suppose that this fluid is in motion and that the flow can be described by a profile of Mach's number $M(y)$ ($M \in ]-1, 1[$). If $\vec{f}$ denotes the source of acoustic waves, we have to solve the transient Galbrun equation in dimensionless form ( $\vec{\xi}$ is the perturbation of lagrangian displacement):
\begin{eqnarray}
\label{equation de Galbrun dans un conduit non borne}
\frac{{\rm D}^2 \vec{\xi}}{{\rm D}t^2}-\vec{\nabla}(\dive\vec{\xi})&=&\vec{f}\ \ {\rm in}\ \ \Omega\times\mathbb{R}_+^* \\
\vec{\xi}\cdot\vec{n}&=&0 \ \ {\rm on}\ \ \partial\Omega\times\mathbb{R}_+^*
\label{condition limite}
\end{eqnarray}   
where $\Dif{~}{t}=\dif{~}{t}+M(y)\dif{~}{x}$ is the material derivative and $\vec{n}$ the unit outward normal vector to $\partial\Omega$.
We complete the equations (\ref{equation de Galbrun dans un conduit non borne})-(\ref{condition limite}) by adding initial conditions.
\section{The Lagrangian vorticity equation}
\noindent When $\vec{f} \in \ahr$, applying the curl operator to (\ref{equation de Galbrun dans un conduit non borne}), we prove that the $\rotscal\vec{\xi}=\psi$ satisfies:
\begin{equation}
\label{equation pour le rotationnel}
 \Diff{\psi}{t}= 2 M'(y) \Dif{~}{t}(\dif{\xi_x}{x})
+\rotscal\vec{f}\ \ {\rm in} \ \ \Omega \times \R_+^* 
\end{equation}
Note that when the flow is uniform ($M(y)$ is constant),
the vorticity $\psi$ can be computed independently of $\vec{\xi}$ and the solution 
is given by~:
\begin{eqnarray*}
\psi(x,\,y,\,t)&= \alpha(x-Mt,\,y)+x\beta(x-Mt,\,y)+\\
  &\frac{1}{M^2}\int_0^x (x-a)(\rotscal\vec{f})(a,\, y,\, t-\frac{x-a}{M})da
\end{eqnarray*}
where $\alpha$ et $\beta$ are two functions that depend only of initial conditions of the problem.
\section{Regularized Formulation Galbrun equation}   
\noindent The idea of the regularization, initially introduced for
Maxwell's equations, was extended to the  time harmonic Galbrun equation by A. S Bonnet-Bendhia and al (2001). The idea  consists in adding the artificial term $s \rot (\rotscal-\psi)$
to Galbrun equation. 

\noindent We replace the initial value problem (\ref{equation de Galbrun dans un conduit non borne})-(\ref{condition limite}) by the equivalent {\it regularized } problem~:
\begin{equation}
{\footnotesize
\label{system de Galbrun regularise}
\left \{\begin{array}{llll}
\Diff{\vec{\xi}}{t}-\vec{\nabla}(\dive\vec{\xi})+s\rot(\rotscal\vec{\xi}-\psi )=\vec{f}&{\rm in}&\Omega\times\mathbb{R_+} \\[0.3cm]
\Dif{~}{t}\left(\Dif{\psi}{t}-2 M'(y)\dif{\xi_x}{x}\right)=\rotscal\vec{f} &{\rm in}&\Omega \times \R_+^*\\[0.3cm]
\vec{\xi}\cdot\vec{n}=0,\ \ \rots\vec{\xi}=\psi&{\rm on}&\  \partial \Omega\times{\mathbb R}_+ 
\end{array}\right.}
\end{equation}
where $s$ is a non negative parameter.\\ \\
In this paper we restrict to the case of uniform flow ($M'(y)=0 \ \
\forall y \in ]-h,\,h[$) for which the problem in $\psi$ and  $\xi_x$
are decoupled. The study of the coupled system will be a subject of future works
\section{Mathematical analysis of the regularized problem}
\noindent When the Mach number is constant the vorticity $\psi$ is known and the problem in  $\vec{\xi}$ is written~: 
\begin{equation}
{\footnotesize
\label{equation de Galbrun regularisee}
\left\{\begin{array}{llll}
\Diff{\vec{\xi}}{t}-\vec{\nabla}(\dive\vec{\xi})+s\rot(\rotscal\vec{\xi})=\vec{f_s}
&{\rm in}&\hspace{-0.2cm}\Omega\times\mathbb{R_+} \\[0.3cm]
\vec{\xi}\cdot\vec{n}=0,\ \ \rots\vec{\xi}=\psi&{\rm on}& \hspace{-0.2cm}\partial \Omega\times{\mathbb R}_+  
\end{array}\right.}
\end{equation}
\noindent where $\vec{f}_s=\vec{f}+s\rot{\psi}$.

\noindent The second boundary condition $(\rots\vec{\xi}=\psi)$ of system (\ref{equation de Galbrun regularisee}) is necessary for the equivalence with the initial problem (\ref{equation de Galbrun dans un conduit non borne})-(\ref{condition limite}). The property of ellipticity of the spatial operator 
$-\vec{\nabla}(\dive\vec)+ s\rot(\rotscal)$ (which equals $-\vec{\Delta}$ if $s=1$) allows to carry out the mathematical and numerical study of this problem 
in a classical frame.

\noindent We consider $\H:={\bf H}_{0}(\Omega) \times \el^2$ where 
$${\bf H}_{0}(\Omega):=\left\{ \vec{\xi}\in(\ah)^2 /\quad \vec{\xi}\cdot\vec{n}=0, \text{ on } \partial\Omega\right \}.
$$
To apply Hille-Yosida theorem, we introduce a new unknown 
\quad $\vec{\zeta}={\rm D} \vec{\xi} / {\rm D} t$. If we pose  $\u=(\vec{\xi},\,\vec{\zeta})^t$, then we can rewrite (\ref{equation de Galbrun regularisee}) under the following form :%

{\small
\begin{equation}
\label{Hille}
\left \{\begin{array}{ll}
\dis\ddt{\u} + A_s \u = \vec{F_s}\\[8pt]
\u(0)=\u_{0}
\end{array}
\right. 
\end{equation}
\begin{alignat*}{2}
A_s\u&=\begin{pmatrix}
 -{\vec{\zeta}}+M\dis\dif{\vec{\xi}}{x} \\[8pt] -
\vec{\nabla}(\dive\vec{\xi})+s\vec{\rot}(\rotscal\vec{\xi})+M\dis\dif{\vec{\zeta}}{x}
    \end{pmatrix}, \\
\vec{F_s} &= \begin{pmatrix}
  0 \\[8pt]\vec{f_s}
  \end{pmatrix}.
\end{alignat*}
}
The domain of the unbounded operator  $A_s$ is defined by~:
$$
\dis D(A_s)=\left\{\u=(\vec{\xi},\,\vec{\zeta})^t \in \H \; \text{tel que}\ \ A_s\u \in \H\right \}.
$$
Using Costabel's identity \cite{Costabel:1991} and Hille-Yosida's theory, we prove that the operator $A_s$  is maximal monotone 
\begin{theorem}:
If $\min(1,\, s)> M^2 $, then for $\vec{f_s} \in C^1(\R_+;\,L^2(\Omega)^2)$ and sufficiently regular initial data, problem 
{\rm (\ref{equation de Galbrun regularisee})} has a unique solution which satisfies:
$$
\vec{\xi}\in C^1(\R_+;\,{\bf H}_0) \cap C^2(\R_+;\,L^2(\Omega)^2)
$$
\end{theorem}
\section{The absorbing boundary conditions}
\noindent For solving numerically the problem, we need to truncate the unbounded domain $\Omega$. The truncated domain $\Omega_b:=]-R;\, R[\times[-h;\, h[ $ requires the introduction of absorbing boundary conditions (ABCs) on the artificial boundaries $\Gamma^{\pm}:=\left \{(x,y) \in \Omega,\,x = \pm R \right \} (R > 0)$. The difficulty is to find appropriate boundary conditions adapted for the regularized formulation. For $s=1$, the boundary conditions that we propose are the following : 
\begin{equation}
  \Dif{\vec{\xi}}{t} + \dif{\vec{\xi}}{n}=\vec{0}\quad  \mbox{ on }  \Gamma^{\pm}
\label{ABC} 
\end{equation}

\noindent where $\vec{n}=(\pm 1,\ 0)$ is the unit outward normal vector of $\Gamma^{\pm}$ and $\vec{\tau}=(0,\ \pm 1)$. The main properties of conditions (\ref{ABC}) are
\begin{itemize}
\item[]i) These are exact conditions for $y$-independent solutions (plane waves) : this is why we speak of first order conditions. 

\item[]ii) Well-posedness: the truncated problem is well posed and for $\vec{f}=\vec{0}$ we have the energy decay result:
$$
\ddt{}E(t)+ \int_{_{\Gamma^{-}\cup \Gamma^{+} }}\left |\dif{\vec{\xi}}{t} \right | ^2 =0 
$$
\noindent where
$$
E(t)=\demi\int_{\Omega_b} \left | \dif{\vec{\xi}}{t} \right | ^2 +
\left |\grad{\vec{\xi}}\right | ^2 - M^2 \left | \dif{\vec{\xi}}{x} \right | ^2 
$$

\item[]iii) These are compatible with a variational formulation of the regularized problem, namely:
\begin{equation}
\label{formulation}
\left. \begin{array}{ll}
\dis\ddtt{~}(\vec{\xi},\,\vec{\eta})+\ddt{}\left ( b(\vec{\xi},\,\vec{\eta}) +c^{\Gamma^{\pm}}(\vec{\xi},\,\vec{\eta})\right )\\[0.3cm]+  
a(\vec{\xi},\,\vec{\eta})+d^{\Gamma^{\pm}}(\vec{\xi},\,\vec{\eta})=(\vec{f},\vec{\eta}) 
\end{array}
\right.
\end{equation}
with $(\cdot,\,\cdot)$ is the $L^2$-inner product and 
{\footnotesize
\begin{alignat*}{2} a(\vec{\xi},\vec{\eta})=&\int_{\Omega_b}\dive\vec{\xi}\dive{\vec{\eta}}+\rotscal\vec{\xi}\rotscal{\vec{\eta}}-M^2\frac{\partial\vec{\xi}}{\partial x}\cdot\frac{\partial{\vec{\eta}}}{\partial x},\ \\
b(\vec{\xi},\vec{\eta})=&\int_{\Omega_b}2M\vec{\xi}\cdot\dif{\vec{\eta}}{x},\, c^{\Gamma^{\pm}(}(\vec{\xi},\vec{\eta})=\int_{\cup\Gamma^{\mp}}(1\pm M)\vec{\xi}\cdot \vec{\eta}\,\mathrm{d}\gamma,\\
d^{\Gamma^{\pm}}(\vec{\xi},\vec{\eta})&=\int_{\Gamma^{-}\cup\Gamma^{+} }{\cal R}\dif{\vec{\xi}}{\vec{\tau }} \cdot \vec{\eta}\,\mathrm{d}\gamma,\ \ 
{\cal R}= 
\left ( \begin{array}{ll}
0 &-1 \\
1&0 
\end{array} 
\right )
\end{alignat*} 
}
\end{itemize}

\begin{remark}
The construction of good absorbing boundary conditions is not trivial. For example the following natural ABC's :

{\small
\begin{eqnarray}
\left .(1-M^2)\,\frac{\partial\xi_1}{\partial x}+\dif{\xi_2}{y}\pm(1\mp M) 
\dif{\xi_1}{t}\;\right|_{\Gamma^{\pm}}&=0
 \\
 \left .(1-M^2)\,\frac{\partial\xi_2}{\partial x}-\dif{\xi_1}{y} \pm(1
   \mp M)\dif{\xi_2}{t}\;\right |_{ \Gamma^{\pm}}&= 0
\end{eqnarray}
}
 are still exact for plane waves and variational, but not stable.  
\end{remark}
\section{Discretization of Galbrun equations}
\noindent 
\noindent The Lagrange finite element approximation of (\ref{formulation}) leads to the following ordinary differential system :
\begin{eqnarray}
\label{semi}
\mathbb{M}_h\frac{{{\rm d}^2\vec{\xi}_h}}{{\rm d} t^2}+( (\mathbb{B}_h + \mathbb{C}^{\Gamma^{\pm}}_h )
\ddt{\vec{\xi}_h}+(\mathbb{A}_h + \mathbb{D}^{\Gamma^{\pm}}_h) \vec{\xi}_h =\vec{F}_h
\end{eqnarray}

\noindent where $\mathbb{M}_h$ is the mass matrix, and $\mathbb{A}_h$, $\mathbb{B}_h$, $\mathbb{C}^{\Gamma^{\pm}}_h$ and $\mathbb{D}^{\Gamma^{\pm}}_h$ are the matrices respectively associated to the bilinear forms $a(.,\,.)$, $b(.,\,.)$, $c^{\Gamma^{\pm}}(.,\,.)$ and $d^{\Gamma^{\pm}}(.,\,.)$.\\ \\
For the time discretization of (\ref{semi}) we use a centered second order finite difference scheme :
\begin{eqnarray*}
\label{shema}
\mathbb{M}_h\frac{\vec{\xi}_h^{n+1} -  2\vec{\xi}_h^{n} + \vec{\xi}_h^{n-1}}{\Delta t ^2}+(\mathbb{B}_h + \mathbb{C}^{\Gamma^{\pm}}_h ) \frac{\vec{\xi}_h^{n+1} - \vec{\xi}_h^{n-1}}{2 \Delta t }+ \\  (\mathbb{A}_h + \mathbb{D}^{\Gamma^{\pm}}_h) \vec{\xi}_h^n =\vec{F}_h^n
\end{eqnarray*}

\section{Numerical Simulation}
\noindent In this first experiment, we simulate a wave initially excited by a rotational source located in the center of the domain, in the presence of a horizontal uniform flow with $M=0.5$.
The first simulation shows that the method is not stable if the equation is not regularized ($s=0$). 
The second simulation corresponds to the regularized case ($s=1$). When can distinguish the two parts of the Lagrangian displacement : the irrotational part corresponds to the outer circular wavefront (whose radius increases with time) while the rotational part corresponds to the inner circular wavefront (whose radius remains constant). Both are centered at a point which is convected by the flow.
\begin{figure}[H]
\centering
  \includegraphics[scale=0.3]{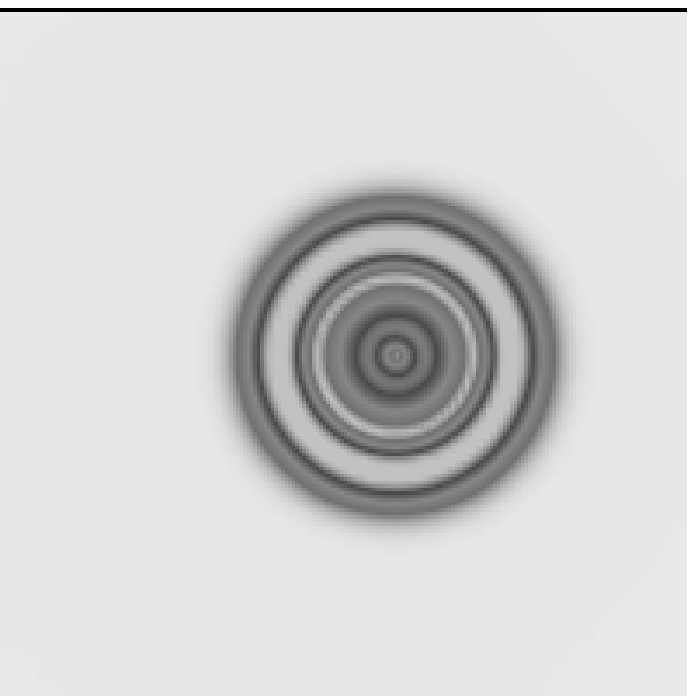}\hspace{0.5cm}
  \includegraphics[scale=0.3]{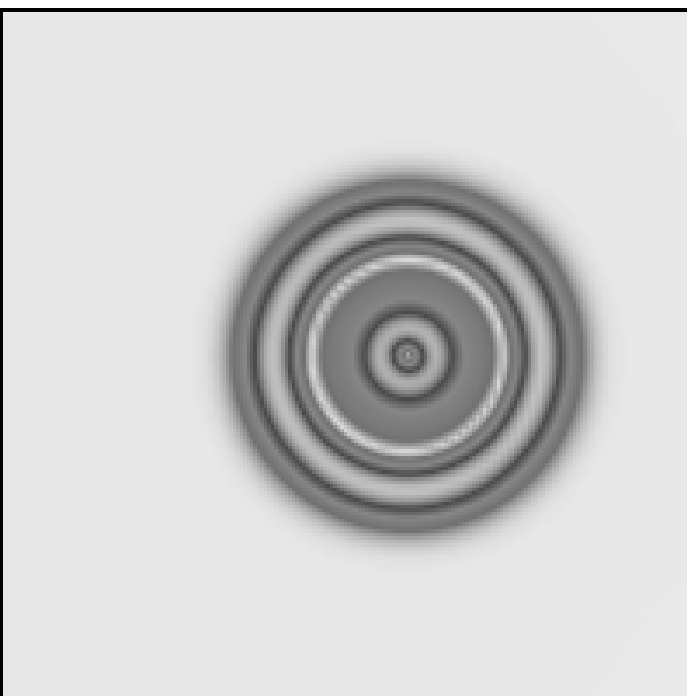}\hspace{0.5cm} 
  \includegraphics[scale=0.3]{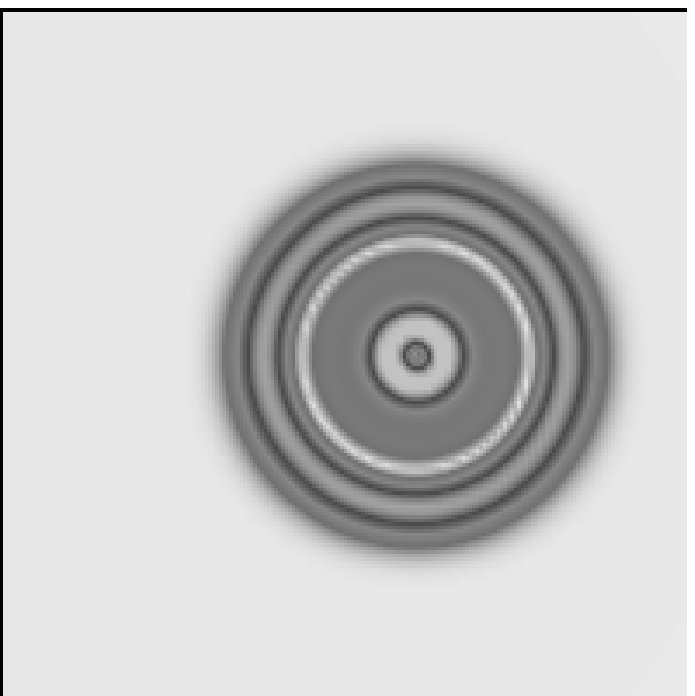}
    \caption{Euclidean norm of $\vec{\xi}$ at $t_1=1.5$s and $t_2=1.75$s and $t_3=2$s. Case $s=1$}
  \label{fig1}
\end{figure}
\noindent In the second experiment, we simulate Galbrun wave propagation in an infinite rigid duct. The wave is excited by a Gaussian signal in time and a quasi-punctual irrotationnal source in space in the uniform flow with $M=0.5$. We remark on this result that the ABC prevents (partially) the unphysical reflexion.
\begin{figure}[H]
\begin{center}
  \includegraphics[scale=0.4]{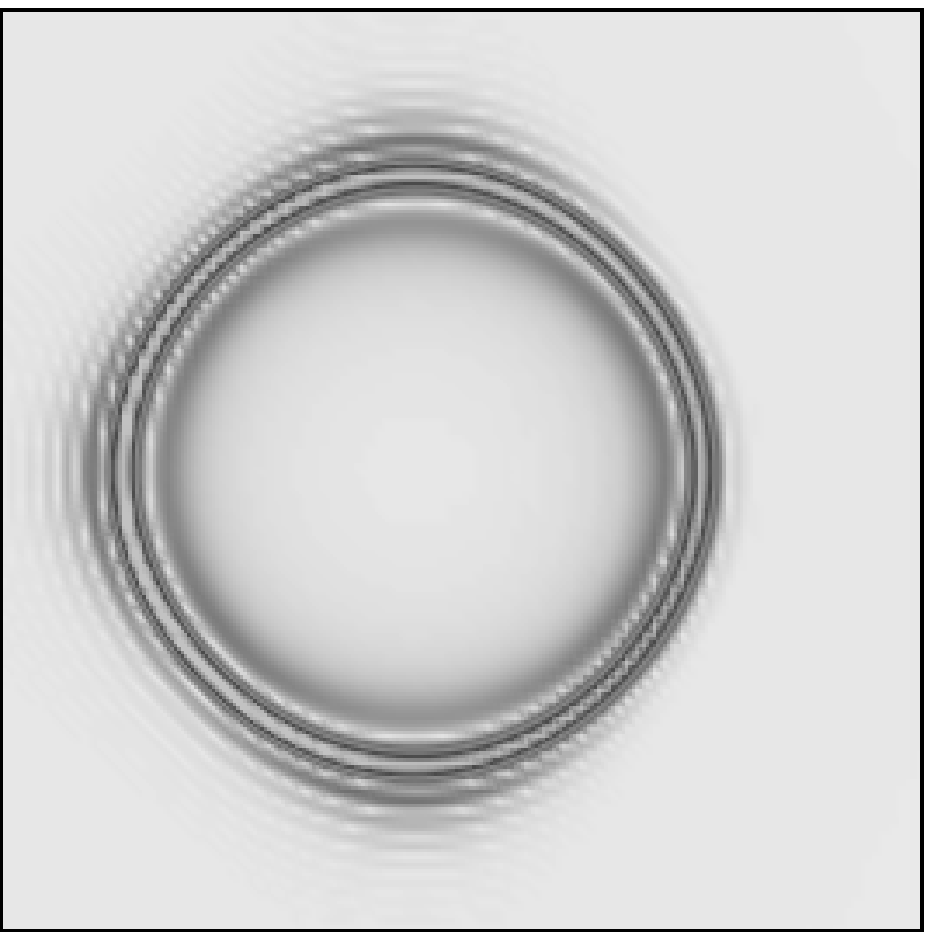}\hspace{0.5cm}
  \includegraphics[scale=0.4]{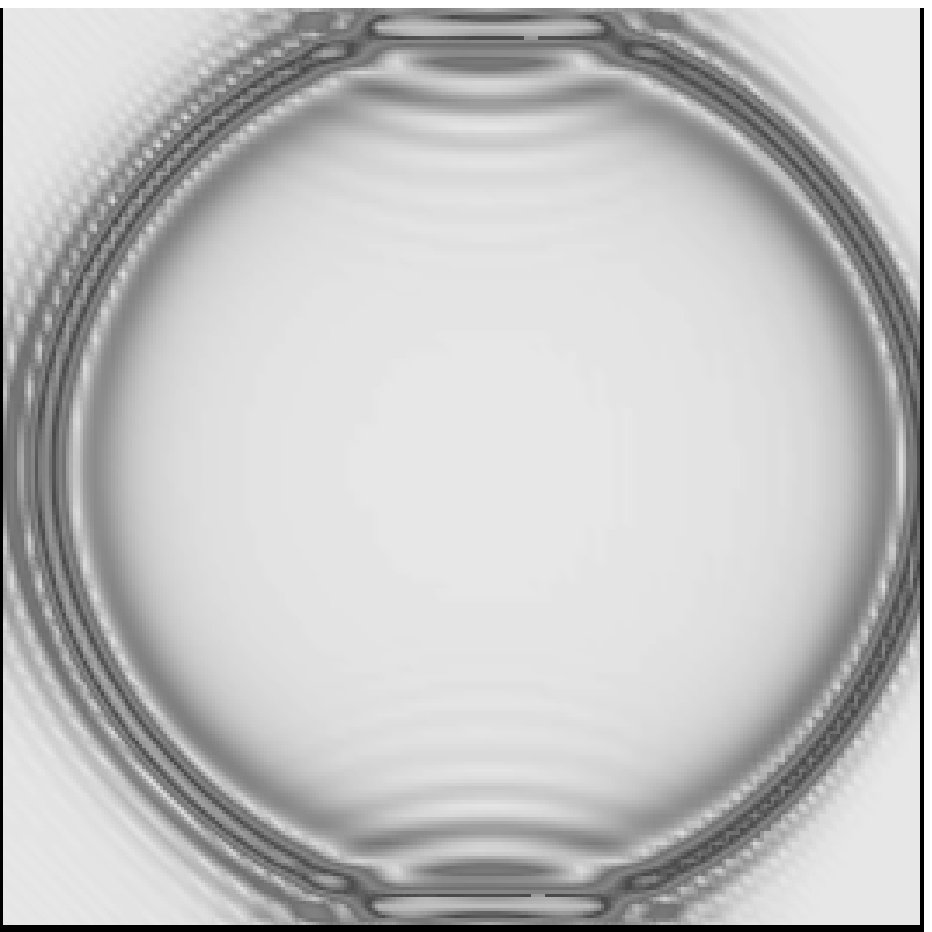}\vspace{0.5cm}
  
  \includegraphics[scale=0.4]{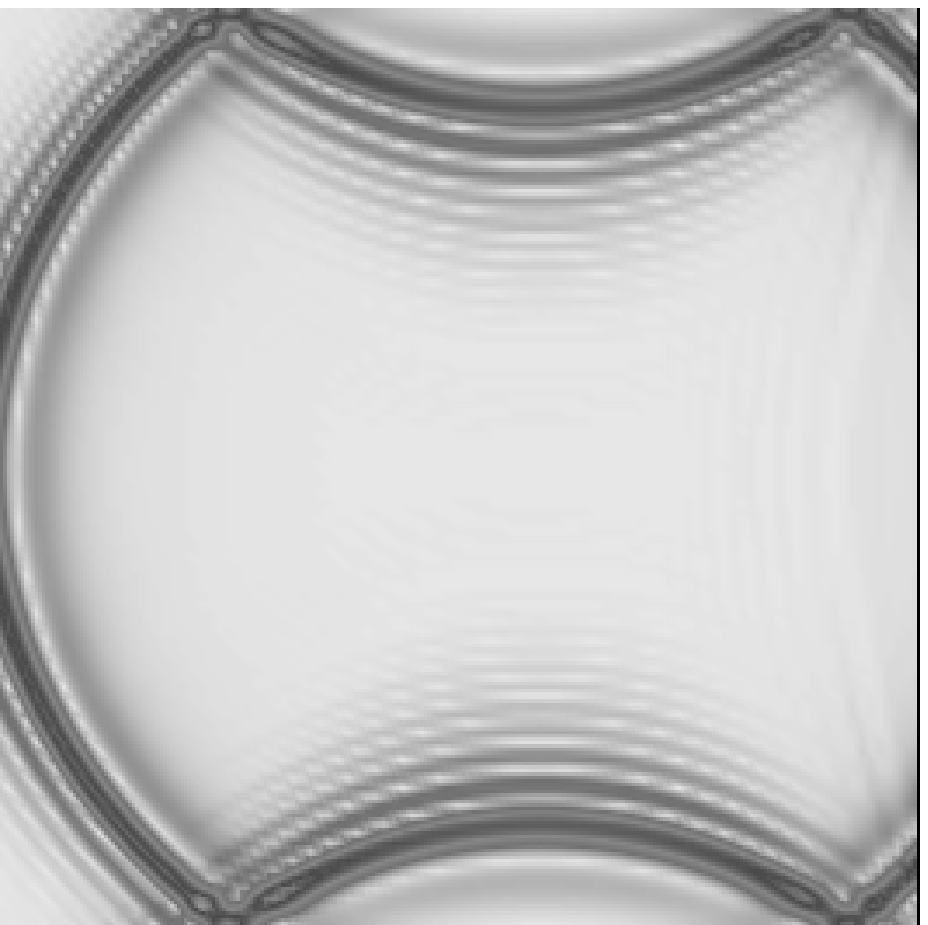}\hspace{0.5cm}
  \includegraphics[scale=0.4]{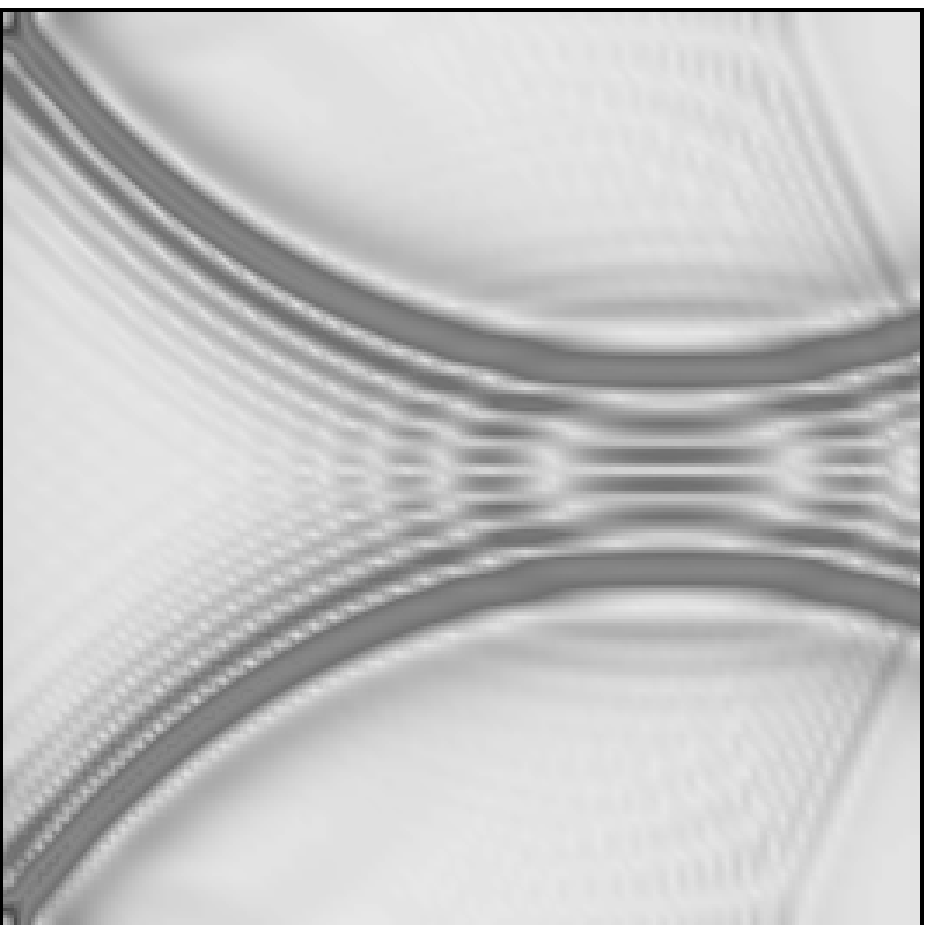}
\end{center}
\caption{Evolution of the Euclidean norm of the displacement $\vec{\xi}$}
\label{cla}
\end{figure}

\end{document}